\documentclass[12pt,a4paper]{article}
\usepackage{indentfirst}
\usepackage{amsfonts}
\usepackage{amssymb}
\usepackage{mathrsfs}
\usepackage{amsmath}
\usepackage{amsthm}
\usepackage{enumerate}
\usepackage{cite}
\usepackage[all]{xy}
\allowdisplaybreaks
\newtheorem{defn}{Definition}[section]
\newtheorem{thm}[defn]{Theorem}
\newtheorem{lem}[defn]{Lemma}
\newtheorem{prop}[defn]{Proposition}
\newtheorem{cor}[defn]{Corollary}
\newtheorem{eg}[defn]{Example}
\newtheorem{re}[defn]{Remark}
\newcommand{\bdefn}{\begin{defn}}
\newcommand{\edefn}{\end{defn}}
\newcommand{\bthm}{\begin{thm}}
\newcommand{\ethm}{\end{thm}}
\newcommand{\blem}{\begin{lem}}
\newcommand{\elem}{\end{lem}}
\newcommand{\bprop}{\begin{prop}}
\newcommand{\eprop}{\end{prop}}
\newcommand{\bcor}{\begin{cor}}
\newcommand{\ecor}{\end{cor}}
\newcommand{\beg}{\begin{eg}}
\newcommand{\eeg}{\end{eg}}
\newcommand{\bre}{\begin{re}}
\newcommand{\ere}{\end{re}}
\newcommand{\bpf}{\begin{proof}}
\newcommand{\epf}{\end{proof}}

\newcommand{\K}{{\rm{\bf K}}}
\newcommand{\id}{{\rm id}}

\newcommand{\Der}{{\rm Der}}

\newcommand{\benu}{\begin{enumerate}}
\newcommand{\eenu}{\end{enumerate}}
\newcommand{\bc}{\begin{center}}
\newcommand{\ec}{\end{center}}
\newcommand{\bea}{\begin{eqnarray}}
\newcommand{\eea}{\end{eqnarray}}
\newcommand{\Bea}{\begin{eqnarray*}}
\newcommand{\Eea}{\end{eqnarray*}}
\newcommand{\beq}{\begin{equation}}
\newcommand{\eeq}{\end{equation}}
\newcommand{\Beq}{\begin{equation*}}
\newcommand{\Eeq}{\end{equation*}}
\newcommand{\bspl}{\begin{split}}
\newcommand{\espl}{\end{split}}

\numberwithin{equation}{section}

\begin{document}
\title{\bf  One-parameter formal deformations of Hom-Lie-Yamaguti algebras}
\author{\normalsize \bf Yao Ma$^1$, Liangyun Chen$^1$,  Jie Lin$^2$}
\date{\small{ $^1$ School of Mathematics and Statistics, Northeast Normal University, Changchun,  130024,  CHINA
  \\$^2$ Sino-European Institute of Aviation Engineering,Civil Aviation University of China,  Tianjin, 300300,  CHINA}}
\maketitle
\begin{abstract}
This paper studies one-parameter formal deformations of Hom-Lie-Yamaguti algebras. The first, second and third cohomology groups on Hom-Lie-Yamaguti algebras extending ones on Lie-Yamaguti algebras are provided. It is proved that first and second cohomology groups are suitable to the deformation theory involving infinitesimals, equivalent deformations and rigidity. However, the third cohomology group is not suitable for the obstructions.
\bigskip

\end{abstract}

\footnote[0]{Corresponding author(L. Chen): chenly640@nenu.edu.cn.}
\footnote[0]{Supported by  NNSF of China (No.11171055 and No.11471090),  NSF of  Jilin province (No.201115006), Scientific Research Foundation for Returned Scholars Ministry of Education of China  and the Fundamental Research Funds for the Central Universities (No.12SSXT139). }

\section{Introduction}
Lie-Yamaguti algebras were introduced by Yamaguti in \cite{Yamaguti2} to give an algebraic interpretation of the characteristic properties of the torsion and curvature of homogeneous spaces with canonical connection in \cite{Nomizu}. He called them ``generalized Lie triple systems'' at first, which were later called ``Lie triple algebras''. Recently, they were renamed as ``Lie-Yamaguti algebras''.

A Hom-type algebra is a kind of algebras whose identities defining the structures are twisted by a linear homomorphism (the twisting map). When the twisting map is the identity map, one recovers the original algebra. The notion of Hom-type algebras was initially introduced in \cite{Hartwig&Larsson&Silvestrov} as ``Hom-Lie algebras'' to describe the $q$-deformation of the Witt and the Virasoro algebras. For more information on various Hom-type algebras one may refer to \cite{Ammar&Mabrouk&Makhlouf, Benayadi&Makhlouf, Hu, Larsson&Silvestrov, Liu&Chen&Ma, Sheng, Sheng&Chen, Yau1, Yau2}. In particular, the notion of Hom-Lie Yamaguti algebras was introduced by Gaparayi and Nourou in \cite{Gaparayi&Nourou1}.

A \textbf{Hom-Lie-Yamaguti algebra} (HLYA for short) is a quadruple $(L, [\cdot,\cdot], \{\cdot, \cdot,\cdot\}, \alpha)$ in which $L$ is a vector space over a field $\K$, ``$[\cdot,\cdot]$'' a binary operation and ``$\{\cdot, \cdot,\cdot\}$'' a ternary operation on $L$, and $\alpha:L\rightarrow L$ a linear map such that for all $x,y,z,u,v\in L$,
\begin{gather}
\alpha([xy])=[\alpha(x)\alpha(y)],\\
\alpha(\{xyz\})=\{\alpha(x)\alpha(y)\alpha(z)\},\\
[xx]=0,\\
\{xxy\}=0,\\
\circlearrowleft_{x,y,z}([[xy]\alpha(z)]+\{xyz\})=0,\label{1.5}\\
\circlearrowleft_{x,y,z}(\{[xy]\alpha(z)\alpha(u)\})=0,\label{1.6}\\
\{\alpha(x)\alpha(y)[uv]\}=[\{xyu\}\alpha^2(v)]+[\alpha^2(u)\{xyv\}],\label{1.7}\\
\{\alpha^2\!(u)\alpha^2\!(v)\{xyz\}\}\!=\!\{\{uvx\}\alpha^2\!(y)\alpha^2\!(z)\}\! +\!\{\alpha^2\!(x)\{uvy\}\alpha^2\!(z)\}\!+\! \{\alpha^2\!(x)\alpha^2\!(y)\{uvz\}\},\label{1.8}
\end{gather}
where $\circlearrowleft_{x,y,z}$ denotes the sum over cyclic permutations of $x,y,z$.

HLYAs generalize Hom-Lie triple systems and Hom-Lie algebras in the same way as Lie-Yamaguti algebras generalize Lie triple systems and Lie algebras, i.e., if $\{xyz\}=0$ for all $x, y, z\in L$, then $(L, [\cdot,\cdot], \alpha)$ becomes a Hom-Lie algebra; if $[xy]=0$ for all $x, y\in L$, then $(L, \{\cdot, \cdot,\cdot\}, \alpha^2)$ becomes a Hom-Lie triple system. When $\alpha=\id_L$, $(L, [\cdot,\cdot], \{\cdot, \cdot,\cdot\})$ becomes a Lie-Yamaguti algebra. There are more examples and properties about HLYA in \cite{Gaparayi&Nourou1, Gaparayi&Nourou2}.

A deformation is a tool to study a mathematical object by deforming it into a family of the same kind of objects depending on a certain parameter. Deformation problems appear in various areas of mathematics, especially in algebra, algebraic and analytic geometry, and mathematical physics. The deformation theory was introduced by Kodaira and Spencer to study complex structures of higher dimensional manifolds (see \cite{Kodaira&Spencer}), which was extended to rings and algebras by Gerstenhaber in \cite{Gerstenhaber1, Gerstenhaber2, Gerstenhaber3, Gerstenhaber4} and to Lie algebras by Nijenhuis and Richardson in \cite{Nijenhuis&Richardson}. They connected deformation theory for associative algebras and Lie algebras with Hochschild cohomology and Chevally-Eilenberg cohomology, respectively. See also \cite{Ammar&Ejbehi&Makhlouf, Elhamdadi&Makhlouf, Flato&Gerstenhaber&Voronov, Kubo&Taniguchi, Lin&Chen&Ma, Ma&Chen&Lin, Makhlouf&Silvestrov} for more deformation theory.

The aim of this paper is to consider the cohomology theory and the one-parameter formal deformation theory of HLYAs based on some work in \cite{Kubo&Taniguchi, Lin&Chen&Ma, Sheng, Yamaguti1}. The rest of this paper is organized as follows. In section 2, we define the first, second and third cohomology groups on HLYAs and show that the first cohomology group corresponds to the derivations space of a HLYA. Section 3 concerns the one-parameter formal deformation theory of HLYAs. We show that the first and second cohomology groups defined in Section 2 fits this one-parameter formal deformation theory but the third one does not.

Throughout this paper $\K$ denotes an arbitrary field.

\section{First, second and third cohomology groups of a Hom-Lie-Yamaguti algebra}
Inspired by the cohomology theory of Lie-Yamaguti algebras in \cite{Yamaguti1}, we introduce the first, second and third cohomology groups of HLYAs.
\bdefn
Let $(L, [\cdot,\cdot], \{\cdot, \cdot,\cdot\}, \alpha)$ be a HLYA. An $n$-linear map $f: \underbrace{L\times \cdots \times L}_{n ~times}\rightarrow L$ is called an $n$-Hom-cochain, if $f$ satisfies
\begin{align}
f(x_{1}, \cdots, x_{2i-1}, x_{2i}, \cdots, x_{n})&=0, ~\text{for}~ x_{2i-1}=x_{2i},\label{cochain1}\\
f(\alpha(x_1), \cdots, \alpha(x_n))&=\alpha\circ f(x_{1}, \cdots, x_{n}).\label{cochain2}
\end{align}
The set of $n$-Hom-cochains is denoted by $HomC^n(L,L)$, for $n\geq1$.
\edefn

\bdefn
Suppose that $(L, [\cdot,\cdot], \{\cdot, \cdot,\cdot\}, \alpha)$ is a HLYA.
\benu[(i)]
\item
A 1-coboundary operator of $(L, [\cdot,\cdot], \{\cdot, \cdot,\cdot\}, \alpha)$ is a pair of maps
\Bea
(\delta_I^1, \delta_{II}^1): HomC^1(L,L)\times HomC^1(L,L)&\longrightarrow& HomC^2(L,L)\times HomC^3(L,L)\\
(f,f)&\longmapsto& (\delta_I^1 f, \delta_{II}^1 f)
\Eea
for $f\in HomC^1(L,L)$ and
\begin{align*}
\delta_I^1 f(x, y)&=[xf(y)]+[f(x)y]-f([xy]),\\
\delta_{II}^1 f(x, y, z)&=\{f(x)yz\}+\{xf(y)z\}+\{xyf(z)\}-f(\{xyz\}).
\end{align*}

\item
A 2-coboundary operator of $(L, [\cdot,\cdot], \{\cdot, \cdot,\cdot\}, \alpha)$ is a pair of maps
\Bea
(\delta_I^2, \delta_{II}^2): HomC^2(L,L)\times HomC^3(L,L)&\longrightarrow& HomC^4(L,L)\times HomC^5(L,L)\\
(f,g)&\longmapsto& (\delta_I^2 f, \delta_{II}^2 g)
\Eea
for $f\in HomC^2(L,L), g\in HomC^3(L,L)$ and
\begin{align*}
\delta_I^2 f(x, y, z, u)=&\{\alpha(x)\alpha(y)f(z, u)\}-f(\{xyz\}, \alpha^2(u))-f(\alpha^2(z), \{xyu\})\\
                         &+g(\alpha(x), \alpha(y), [zu])-[\alpha^2(z)g(x, y, u)]-[g(x, y, z)\alpha^2(u)],\\
\delta_{II}^2 g(x, y, u, v, w)=&\{\alpha^2(x)\alpha^2(y)g(u, v, w)\}-\{g(x, y, u)\alpha^2(v)\alpha^2(w)\}\\
                               &-\{\alpha^2(u)g(x, y, v)\alpha^2(w)\}-\{\alpha^2(u)\alpha^2(v)g(x, y, w)\}\\
                               &+g(\alpha^2(x), \alpha^2(y), \{uvw\})-g(\{xyu\}, \alpha^2(v), \alpha^2(w))\\
                               &-g(\alpha^2(u), \{xyv\}, \alpha^2(w))-g(\alpha^2(u), \alpha^2(v), \{xyw\}).
\end{align*}

\item
A 3-coboundary operator of $(L, [\cdot,\cdot], \{\cdot, \cdot,\cdot\}, \alpha)$ is a pair of maps
\Bea
(\delta_I^3, \delta_{II}^3): HomC^4(L,L)\times HomC^5(L,L)&\longrightarrow& HomC^6(L,L)\times HomC^7(L,L)\\
(f,g)&\longmapsto& (\delta_I^3 f, \delta_{II}^3 g)
\Eea
for $f\in HomC^4(L,L), g\in HomC^5(L,L)$ and
\begin{align*}
 &\delta_I^3 f(x_1, \cdots, x_6)\\
=&\{\alpha^3(x_1)\alpha^3(x_2)f(x_3, \cdots, x_6)\}-\{\alpha^3(x_3)\alpha^3(x_4)f(x_1, x_2, x_5, x_6)\}\\
 &+\sum_{k=1}^2\sum_{i=2k+1}^6(-1)^k f(\alpha^2(x_1),\cdots, \widehat{\alpha^2(x_{2k-1})}, \widehat{\alpha^2(x_{2k})}, \cdots, \{x_{2k-1}x_{2k}x_{i}\},\cdots, \alpha^2(x_6))\\
 &-g(\alpha(x_1),\cdots, \alpha(x_4), [x_5x_6])+[\alpha^4(x_5)g(x_1,\cdots,x_4,x_6)] +[g(x_1,\cdots,x_5)\alpha^4(x_6)],
\end{align*}
\begin{align*}
&\delta_{II}^3 g(x_1, \cdots, x_7)\\
=&\sum_{k=1}^3(-1)^{k+1}\{\alpha^4(x_{2k-1})\alpha^4(x_{2k})g(x_1, \cdots, \widehat{x_{2k-1}}, \widehat{x_{2k}}, \cdots, x_7)\\
&+\sum_{k=1}^3\sum_{i=2k+1}^7(-1)^k g(\alpha^2(x_1),\cdots, \widehat{\alpha^2(x_{2k-1})}, \widehat{\alpha^2(x_{2k})}, \cdots, \{x_{2k-1}x_{2k}x_{i}\},\cdots, \alpha^2(x_7))\\
&+\{g(x_1,\cdots, x_5)\alpha^4(x_6)\alpha^4(x_7)\}-\{g(x_1,\cdots, x_4, x_6)\alpha^4(x_5)\alpha^4(x_7)\},
\end{align*}
where the sign \textasciicircum ~indicates that the element below must be omitted.
\eenu
\edefn

\bthm
The coboundary operators $(\delta_I^i, \delta_{II}^i)$ are well defined, for $i=1, 2, 3$.
\ethm
\bpf
Take $(f, f)\in HomC^1(L,L)\times HomC^1(L,L)$. It is clear that $\delta_I^1 f$ and $\delta_{II}^1 f$ satisfy (\ref{cochain1}). Notice that
\begin{align*}
\delta_I^1 f(\alpha(x), \alpha(y))=&[\alpha(x)f(\alpha(y))]+[f(\alpha(x))\alpha(y)]-f([\alpha(x)\alpha(y)])\\
                                  =&\alpha([xf(y)])+\alpha([f(x)y])-\alpha(f([xy]))=\alpha\circ\delta_I^1 f(x, y)
\intertext{and}
\delta_{II}^1 f(\alpha(x), \alpha(y), \alpha(z))=&\{f(\alpha(x))\alpha(y)\alpha(z)\}+\{\alpha(x)f(\alpha(y))\alpha(z)\}+\{\alpha(x)\alpha(y)f(\alpha(z))\}\\
                                                 &-f(\{\alpha(x)\alpha(y)\alpha(z)\})\\
                                                =&\alpha(\{f(x)yz\}+\{xf(y)z\}+\{xyf(z)\}-f(\{xyz\}))\\
                                                =&\alpha\circ\delta_{II}^1 f(x, y, z).
\end{align*}
Then $(\delta_I^1, \delta_{II}^1)$ is well defined.

Now let $(f, g)\in HomC^2(L,L)\times HomC^3(L,L)$. Then $\delta_I^2 f$ and $\delta_{II}^2 g$ satisfy (\ref{cochain1}) and
\begin{align*}
 &\delta_I^2 f(\alpha(x), \alpha(y), \alpha(z), \alpha(u))\\
=&\{\alpha^2(x)\alpha^2(y)f(\alpha(z), \alpha(u))\}-f(\alpha(\{xyz\}), \alpha^3(u))-f(\alpha^3(z), \alpha(\{xyu\}))\\
 &+g(\alpha^2(x), \alpha^2(y), \alpha([zu]))-[\alpha^3(z)g(\alpha(x), \alpha(y), \alpha(u))]-[g(\alpha(x), \alpha(y), \alpha(z))\alpha^3(u)]\\
=&\alpha(\{\alpha(x)\alpha(y)f(z, u)\})-\alpha\circ f(\{xyz\}, \alpha^2(u))-\alpha\circ f(\alpha^2(z), \{xyu\})\\
 &+\alpha\circ g(\alpha(x), \alpha(y), [zu])-\alpha([\alpha^2(z)g(x, y, u)])-\alpha([g(x, y, z)\alpha^2(u)])\\
=&\alpha\circ \delta_I^2 f(x, y, z, u).
\end{align*}
Similarly,
$$\delta_{II}^2 g(\alpha(x), \alpha(y), \alpha(u), \alpha(v), \alpha(w))=\alpha\circ \delta_{II}^2 g(x, y, u, v, w).$$

One proves $(\delta_I^3 f, \delta_{II}^3 g)\in HomC^6(L,L)\times HomC^7(L,L)$ if $(f, g)\in HomC^4(L,L)\times HomC^5(L,L)$ in the same way. Hence the theorem follows.
\epf

Moreover, for $(f, g)\in HomC^2(L,L)\times HomC^3(L,L)$, we define another 2-coboundary operator of $(L, [\cdot,\cdot], \{\cdot, \cdot,\cdot\}, \alpha)$ as
\Bea
(d_I^2, d_{II}^2): HomC^2(L,L)\times HomC^3(L,L)&\longrightarrow& HomC^3(L,L)\times HomC^4(L,L)\\
(f,g)&\longmapsto& (d_I^2 f, d_{II}^2 g)
\Eea
where
\begin{align*}
d_I^2 f(x, y, z)=&\circlearrowleft_{x,y,z}([f(x, y)\alpha(z)]+f([xy], \alpha(z))+g(x, y, z)),\\
d_{II}^2 g(x, y, z, u)=&\circlearrowleft_{x,y,z}(\{f(x, y)\alpha(z)\alpha(u)\}+g([xy], \alpha(z), \alpha(u))).
\end{align*}
It is easy to prove that $(d_I^2, d_{II}^2)$ is well defined.

\bthm\label{delta^2=0}
With notations as above, we have
$$(\delta_I^2, \delta_{II}^2)(\delta_I^1, \delta_{II}^1)=(0,0), \quad (d_I^2, d_{II}^2)(\delta_I^1, \delta_{II}^1)=(0,0), \quad (\delta_I^3, \delta_{II}^3)(\delta_I^2, \delta_{II}^2)=0.$$
\ethm
\bpf
Suppose $(f, f)\in HomC^1(L,L)\times HomC^1(L,L)$. Then
$$(\delta_I^2, \delta_{II}^2)(\delta_I^1, \delta_{II}^1)(f, f)=(\delta_I^2, \delta_{II}^2)(\delta_I^1 f, \delta_{II}^1 f)=(\delta_I^2\delta_I^1 f, \delta_{II}^2\delta_{II}^1 f).$$
Using (\ref{1.7}), (\ref{1.8}) and (\ref{cochain2}), we have
\begin{align*}
&\delta_I^2\delta_I^1 f(x, y, z, u)\\
=&\{\alpha(x)\alpha(y)\delta_I^1 f(z, u)\}-\delta_I^1 f(\{xyz\}, \alpha^2(u))-\delta_I^1 f(\alpha^2(z), \{xyu\})\\
 &+\delta_{II}^1 f(\alpha(x), \alpha(y), [zu])-[\alpha^2(z)\delta_{II}^1 f(x, y, u)]-[\delta_{II}^1 f(x, y, z)\alpha^2(u)]\\
=&\{\alpha(x)\alpha(y)[zf(u)]\}+\{\alpha(x)\alpha(y)[f(z)u]\}-\{\alpha(x)\alpha(y)f([zu])\}\\
 &-[\{xyz\}f\alpha^2(u)]-[f(\{xyz\})\alpha^2(u)]+f([\{xyz\}\alpha^2(u)])\\
 &-[\alpha^2(z)f(\{xyu\})]-[f\alpha^2(z)\{xyu\}]+f([\alpha^2(z)\{xyu\}])\\
 &+\{f\alpha(x)\alpha(y)[zu]\}+\{\alpha(x)f\alpha(y)[zu]\}+\{\alpha(x)\alpha(y)f([zu])\}-f(\{\alpha(x)\alpha(y)[zu]\})\\
 &-[\alpha^2(z)\{f(x)yu\}]-[\alpha^2(z)\{xf(y)u\}]-[\alpha^2(z)\{xyf(u)\}]+[\alpha^2(z)f(\{xyu\})]\\
 &-[\{f(x)yz\}\alpha^2(u)]-[\{xf(y)z\}\alpha^2(u)]-[\{xyf(z)\}\alpha^2(u)]+[f(\{xyz\})\alpha^2(u)]\\
=&0
\end{align*}
and
\begin{align*}
&\delta_{II}^2\delta_{II}^1 f(x, y, u, v, w)\\
=&\{\alpha^2(x)\alpha^2(y)\delta_{II}^1 f(u, v, w)\}-\{\delta_{II}^1 f(x, y, u)\alpha^2(v)\alpha^2(w)\}-\{\alpha^2(u)\delta_{II}^1 f(x, y, v)\alpha^2(w)\}\\
&-\{\alpha^2(u)\alpha^2(v)\delta_{II}^1 f(x, y, w)\}+\delta_{II}^1 f(\alpha^2(x), \alpha^2(y), \{uvw\})-\delta_{II}^1 f(\{xyu\}, \alpha^2(v), \alpha^2(w))\\
&-\delta_{II}^1 f(\alpha^2(u), \{xyv\}, \alpha^2(w))-\delta_{II}^1 f(\alpha^2(u), \alpha^2(v), \{xyw\})\\
=&\{\alpha^2(x)\alpha^2(y)\{f(u)vw\}\}+\{\alpha^2(x)\alpha^2(y)\{uf(v)w\}\}+\{\alpha^2(x)\alpha^2(y)\{uvf(w)\}\}\\
 &-\{\alpha^2(x)\alpha^2(y)f(\{uvw\})\}-\{\{f(x)yu\}\alpha^2(v)\alpha^2(w)\}-\{\{xf(y)u\}\alpha^2(v)\alpha^2(w)\}\\
 &-\{\{xyf(u)\}\alpha^2(v)\alpha^2(w)\}+\{f(\{xyu\})\alpha^2(v)\alpha^2(w)\}-\{\alpha^2(u)\{f(x)yv\}\alpha^2(w)\}\\
 &-\{\alpha^2(u)\{xf(y)v\}\alpha^2(w)\}-\{\alpha^2(u)\{xyf(v)\}\alpha^2(w)\}+\{\alpha^2(u)f(\{xyv\})\alpha^2(w)\}\\
 &-\{\alpha^2(u)\alpha^2(v)\{f(x)yw\}\}-\{\alpha^2(u)\alpha^2(v)\{xf(y)w\}\}-\{\alpha^2(u)\alpha^2(v)\{xyf(w)\}\}\\
 &+\{\alpha^2(u)\alpha^2(v)f(\{xyw\})\}+\{f\alpha^2(x)\alpha^2(y)\{uvw\}\}+\{\alpha^2(x)f\alpha^2(y)\{uvw\}\}\\
 &+\{\alpha^2(x)\alpha^2(y)f(\{uvw\})\}-f(\{\alpha^2(x)\alpha^2(y)\{uvw\}\})-\{f(\{xyu\})\alpha^2(v)\alpha^2(w)\}\\
 &-\{\{xyu\}f\alpha^2(v)\alpha^2(w)\}-\{\{xyu\}\alpha^2(v)f\alpha^2(w)\}+f(\{\{xyu\}\alpha^2(v)\alpha^2(w)\})\\
 &-\{f\alpha^2(u)\{xyv\}\alpha^2(w)\}-\{\alpha^2(u)f(\{xyv\})\alpha^2(w)\}-\{\alpha^2(u)\{xyv\}f\alpha^2(w)\}\\
 &+f(\{\alpha^2(u)\{xyv\}\alpha^2(w)\})-\{f\alpha^2(u)\alpha^2(v)\{xyw\}\}-\{\alpha^2(u)f\alpha^2(v)\{xyw\}\}\\
 &-\{\alpha^2(u)\alpha^2(v)f(\{xyw\})\}+f(\{\alpha^2(u)\alpha^2(v)\{xyw\}\})\\
=&0.
\end{align*}
Moreover, by (\ref{1.5}) and (\ref{1.6}), one obtains
\begin{align*}
d_I^2\delta_I^1 f(x, y, z)
=&\circlearrowleft_{x,y,z}([\delta_I^1 f(x, y)\alpha(z)]+\delta_I^1 f([xy], \alpha(z))+\delta_{II}^1 f(x, y, z))\\
=&\circlearrowleft_{x,y,z}([[xf(y)]\alpha(z)]+[[f(x)y]\alpha(z)]-[f([xy])\alpha(z)])\\
&+\circlearrowleft_{x,y,z}([[xy]f\alpha(z)]+[f([xy])\alpha(z)]-f([[xy]\alpha(z)]))\\
&+\circlearrowleft_{x,y,z}(\{f(x)yz\}+\{xf(y)z\}+\{xyf(z)\}-f(\{xyz\}))\\
=&0
\end{align*}
and
\begin{align*}
&d_{II}^2\delta_{II}^1 f(x, y, z, u)\\
=&\circlearrowleft_{x,y,z}(\{\delta_I^1 f(x, y)\alpha(z)\alpha(u)\}+\delta_{II}^1 f([xy], \alpha(z), \alpha(u)))\\
=&\circlearrowleft_{x,y,z}(\{[xf(y)]\alpha(z)\alpha(u)\}+\{[f(x)y]\alpha(z)\alpha(u)\}-\{f([xy])\alpha(z)\alpha(u)\})\\
 &+\!\circlearrowleft_{x,y,z}(\{f([xy])\alpha(z)\alpha(u)\}+\{[xy]f\alpha(z)\alpha(u)\}+\{[xy]\alpha(z)f\alpha(u)\}\!-\!f(\{[xy]\alpha(z)\alpha(u)\}))\\
=&\{[xf(y)]\alpha(z)\alpha(u)\}+\{[yf(z)]\alpha(x)\alpha(u)\}+\{[zf(x)]\alpha(y)\alpha(u)\}\\
 &+\{[f(x)y]\alpha(z)\alpha(u)\}+\{[f(y)z]\alpha(x)\alpha(u)\}+\{[f(z)x]\alpha(y)\alpha(u)\}\\
 &+\{[xy]f\alpha(z)\alpha(u)\}+\{[yz]f\alpha(x)\alpha(u)\}+\{[zx]f\alpha(y)\alpha(u)\}\\
=&0.
\end{align*}

For $(f, g)\in HomC^2(L,L)\times HomC^3(L,L)$,
\begin{align*}
 &\delta_{I}^3\delta_{I}^2f(x_1, \cdots, x_6)\\
=&\{\alpha^3(x_1)\alpha^3(x_2)\delta_{I}^2f(x_3, \cdots, x_6)\}-\{\alpha^3(x_3)\alpha^3(x_4)\delta_{I}^2f(x_1, x_2, x_5, x_6)\}\\
 +&\sum_{k=1}^2\sum_{i=2k+1}^6(-1)^k \delta_{I}^2f(\alpha^2(x_1),\cdots, \widehat{\alpha^2(x_{2k-1})}, \widehat{\alpha^2(x_{2k})}, \cdots, \{x_{2k-1}x_{2k}x_{i}\},\cdots, \alpha^2(x_6))\\
 -&\delta_{I}^2g(\alpha(x_1),\cdots, \alpha(x_4), [x_5x_6])+[\alpha^4(x_5)\delta_{I}^2g(x_1,\cdots,x_4,x_6)] +[\delta_{I}^2g(x_1,\cdots,x_5)\alpha^4(x_6)]\\
=&\{\alpha^3(x_1)\alpha^3(x_2)\{\alpha(x_3)\alpha(x_4)f(x_5, x_6)\}\}-\{\alpha^3(x_3)\alpha^3(x_4)\{\alpha(x_1)\alpha(x_2)f(x_5, x_6)\}\}\\
 -&\{\{\alpha(x_1)\alpha(x_2)\alpha(x_3)\}\alpha^3(x_4)\alpha^2f(x_5, x_6)\}-\{\alpha^3(x_3)\{\alpha(x_1)\alpha(x_2)\alpha(x_4)\}\alpha^2f(x_5, x_6)\}\\
 +&f(\{\{x_1x_2x_3\}\alpha^2(x_4)\alpha^2(x_5)\}, \alpha^4(x_6))+f(\{\alpha^2(x_3)\{x_1x_2x_4\}\alpha^2(x_5)\}, \alpha^4(x_6))\\
 +&f(\{\alpha^2(x_3)\alpha^2(x_4)\{x_1x_2x_5\}\}, \alpha^4(x_6))-f(\{\alpha^2(x_1)\alpha^2(x_2)\{x_3x_4x_5\}\}, \alpha^4(x_6))\\
 +&f(\alpha^4(x_5), \{\{x_1x_2x_3\}\alpha^2(x_4)\alpha^2(x_6)\})+f(\alpha^4(x_5), \{\alpha^2(x_3)\{x_1x_2x_4\}\alpha^2(x_6)\})\\
 +&f(\alpha^4(x_5), \{\alpha^2(x_3)\alpha^2(x_4)\{x_1x_2x_6\}\})-f(\alpha^4(x_5), \{\alpha^2(x_1)\alpha^2(x_2)\{x_3x_4x_6\}\})\\
 -&\{\alpha^3(x_1)\alpha^3(x_2)[\alpha^2(x_5)g(x_3, x_4, x_6)]\}+[\{\alpha^2(x_1)\alpha^2(x_2)\alpha^2(x_5)\}\alpha^2g(x_3, x_4, x_6)]\\
 +&[\alpha^4(x_5)\{\alpha^2(x_1)\alpha^2(x_2)g(x_3, x_4, x_6)\}]-\{\alpha^3(x_1)\alpha^3(x_2)[g(x_3, x_4, x_5)\alpha^2(x_6)]\}\\
 +&[\alpha^2g(x_3, x_4, x_5)\{\alpha^2(x_1)\alpha^2(x_2)\alpha^2(x_6)\}]+[\{\alpha^2(x_1)\alpha^2(x_2)g(x_3, x_4, x_5)\}\alpha^4(x_6)]\\
 +&\{\alpha^3(x_3)\alpha^3(x_4)[\alpha^2(x_5)g(x_1, x_2, x_6)]\}-[\{\alpha^2(x_3)\alpha^2(x_4)\alpha^2(x_5)\}\alpha^2g(x_1, x_2, x_6)]\\
 -&[\alpha^4(x_5)\{\alpha^2(x_3)\alpha^2(x_4)g(x_1, x_2, x_6)\}]+\{\alpha^3(x_3)\alpha^3(x_4)[g(x_1, x_2, x_5)\alpha^2(x_6)]\}\\
 -&[\alpha^2g(x_1, x_2, x_5)\{\alpha^2(x_3)\alpha^2(x_4)\alpha^2(x_6)\}]-[\{\alpha^2(x_3)\alpha^2(x_4)g(x_1, x_2, x_5)\}\alpha^4(x_6)]\\
 -&g(\alpha^3(x_3), \alpha^3(x_4), [\{x_1x_2x_5\}\alpha^2(x_6)])-g(\alpha^3(x_3), \alpha^3(x_4), [\alpha^2(x_5)\{x_1x_2x_6\}])\\
 +&g(\alpha^3(x_3), \alpha^3(x_4), \{\alpha(x_1)\alpha(x_2)[x_5x_6]\})+g(\alpha^3(x_1), \alpha^3(x_2), [\{x_3x_4x_5\}\alpha^2(x_6)])\\
 +&g(\alpha^3(x_1), \alpha^3(x_2), [\alpha^2(x_5)\{x_3x_4x_6\}])-g(\alpha^3(x_1), \alpha^3(x_2), \{\alpha(x_3)\alpha(x_4)[x_5x_6]\})\\
 +&\{\alpha g(x_1, x_2, x_3)\alpha^3(x_4)[\alpha^2(x_5)\alpha^2(x_6)]\}-[\alpha^4(x_5)\{g(x_1, x_2, x_3)\alpha^2(x_4)\alpha^2(x_6)\}]\\
 -&[\{g(x_1, x_2, x_3)\alpha^2(x_4)\alpha^2(x_5)\}\alpha^4(x_6)]+\{\alpha^3(x_3)\alpha g(x_1, x_2, x_4)[\alpha^2(x_5)\alpha^2(x_6)]\}\\
 -&[\alpha^4(x_5)\{\alpha^2(x_3)g(x_1, x_2, x_4)\alpha^2(x_6)\}]-[\{\alpha^2(x_3)g(x_1, x_2, x_4)\alpha^2(x_5)\}\alpha^4(x_6)]\\
=&0
\intertext{and}
 &\delta_{II}^3\delta_{II}^2g(x_1, \cdots, x_7)\\
=&\sum_{k=1}^3(-1)^{k+1}\{\alpha^4(x_{2k-1})\alpha^4(x_{2k})\delta_{II}^2g(x_1, \cdots, \widehat{x_{2k-1}}, \widehat{x_{2k}}, \cdots, x_7)\\
 +&\sum_{k=1}^3\sum_{i=2k+1}^7(-1)^k \delta_{II}^2g(\alpha^2(x_1),\cdots, \widehat{\alpha^2(x_{2k-1})}, \widehat{\alpha^2(x_{2k})}, \cdots, \{x_{2k-1}x_{2k}x_{i}\},\cdots, \alpha^2(x_7))\\
 +&\{\delta_{II}^2g(x_1,\cdots, x_5)\alpha^4(x_6)\alpha^4(x_7)\}-\{\delta_{II}^2g(x_1,\cdots, x_4, x_6)\alpha^4(x_5)\alpha^4(x_7)\}\\
=&\{\alpha^4(x_1)\alpha^4(x_2)\{\alpha^2(x_3)\alpha^2(x_4)g(x_5,\! x_6,\! x_7)\}\}\!-\!\{\alpha^4(x_3)\alpha^4(x_4)\{\alpha^2(x_1)\alpha^2(x_2)g(x_5,\! x_6,\! x_7)\}\}\\
 -&\{\alpha^2(\{x_1x_2x_3\})\alpha^4(x_4)\alpha^2g(x_5, x_6, x_7)\}-\{\alpha^4(x_3)\alpha^2(\{x_1x_2x_4\})\alpha^2g(x_5, x_6, x_7)\}\\
 -&\{\alpha^4(x_1)\alpha^4(x_2)\{g(x_3, x_4, x_5)\alpha^2(x_6)\alpha^2(x_7)\}\}+\{\alpha^2g(x_3, x_4, x_5)\alpha^2(\{x_1x_2x_6\})\alpha^4(x_7)\}\\
 +&\{\alpha^2g(x_3, x_4, x_5)\alpha^4(x_6)\alpha^2(\{x_1x_2x_7\})\}+\{\{\alpha^2(x_1)\alpha^2(x_2)g(x_3, x_4, x_5)\}\alpha^4(x_6)\alpha^4(x_7)\}\\
 -&\{\alpha^4(x_1)\alpha^4(x_2)\{\alpha^2(x_5)g(x_3, x_4, x_6)\alpha^2(x_7)\}\}+\{\alpha^2(\{x_1x_2x_5\})\alpha^2g(x_3, x_4, x_6)\alpha^4(x_7)\}\\
 +&\{\alpha^4(x_5)\alpha^2g(x_3, x_4, x_6)\alpha^2(\{x_1x_2x_7\})\}+\{\alpha^4(x_5)\{\alpha^2(x_1)\alpha^2(x_2)g(x_3, x_4, x_6)\}\alpha^4(x_7)\}\\
 -&\{\alpha^4(x_1)\alpha^4(x_2)\{\alpha^2(x_5)\alpha^2(x_6)g(x_3,\! x_4,\! x_7)\}\}\!+\!\{\alpha^4(x_5)\alpha^4(x_6)\{\alpha^2(x_1)\alpha^2(x_2)g(x_3,\! x_4,\! x_7)\}\}\\
 +&\{\alpha^2(\{x_1x_2x_5\})\alpha^4(x_6)\alpha^2g(x_3, x_4, x_7)\}+\{\alpha^4(x_5)\alpha^2(\{x_1x_2x_6\})\alpha^2g(x_3, x_4, x_7)\}\\
 +&\{\alpha^4(x_3)\alpha^4(x_4)\{g(x_1, x_2, x_5)\alpha^2(x_6)\alpha^2(x_7)\}\}-\{\alpha^2g(x_1, x_2, x_5)\alpha^2(\{x_3x_4x_6\})\alpha^4(x_7)\}\\
 -&\{\alpha^2g(x_1, x_2, x_5)\alpha^4(x_6)\alpha^2(\{x_3x_4x_7\})\}-\{\{\alpha^2(x_3)\alpha^2(x_4)g(x_1, x_2, x_5)\}\alpha^4(x_6)\alpha^4(x_7)\}\\
 +&\{\alpha^4(x_3)\alpha^4(x_4)\{\alpha^2(x_5)g(x_1, x_2, x_6)\alpha^2(x_7)\}\}-\{\alpha^2(\{x_3x_4x_5\})g(x_1, x_2, x_6)\alpha^4(x_7)\}\\
 -&\{\alpha^4(x_5)\alpha^2g(x_1, x_2, x_6)\alpha^2(\{x_3x_4x_7\})\}-\{\alpha^4(x_5)\{\alpha^2(x_3)\alpha^2(x_4)g(x_1, x_2, x_6)\}\alpha^4(x_7)\}\\
 +&\{\alpha^4(x_3)\alpha^4(x_4)\{\alpha^2(x_5)\alpha^2(x_6)g(x_1,\! x_2,\! x_7)\}\}\!-\!\{\alpha^4(x_5)\alpha^4(x_6)\{\alpha^2(x_3)\alpha^2(x_4)g(x_1,\! x_2,\! x_7)\}\}\\
 -&\{\alpha^2(\{x_3x_4x_5\})\alpha^4(x_6)\alpha^2g(x_1, x_2, x_7)\}-\{\alpha^4(x_5)\alpha^2(\{x_3x_4x_6\})\alpha^2g(x_1, x_2, x_7)\}\\
 -&\{\alpha^4(x_5)\alpha^4(x_6)\{g(x_1, x_2, x_3)\alpha^2(x_4)\alpha^2(x_7)\}\}+\{\alpha^2g(x_1, x_2, x_3)\alpha^4(x_4)\alpha^2(\{x_5x_6x_7\})\}\\
 -&\{\{g(x_1,\! x_2,\! x_3)\alpha^2(x_4)\alpha^2(x_5)\}\alpha^4(x_6)\alpha^4(x_7)\}\!-\!\{\{\alpha^4(x_5)\{g(x_1,\! x_2,\! x_3)\alpha^2(x_4)\alpha^2(x_6)\}\alpha^4(x_7)\}\\
 -&\{\alpha^4(x_5)\alpha^4(x_6)\{\alpha^2(x_3)g(x_1, x_2, x_4)\alpha^2(x_7)\}\}+\{\alpha^4(x_3)\alpha^2g(x_1, x_2, x_4)\alpha^2(\{x_5x_6x_7\})\}\\
 -&\{\{\alpha^2(x_3)g(x_1,\! x_2,\! x_4)\alpha^2(x_5)\}\alpha^4(x_6)\alpha^4(x_7)\}\!-\! \{\alpha^4(x_5)\{\alpha^2(x_3)g(x_1,\! x_2,\! x_4)\alpha^2(x_6)\}\alpha^4(x_7)\}\\
 +&g(\{\{x_1x_2x_3\}\alpha^2(x_4)\alpha^2(x_5)\},\! \alpha^4(x_6),\! \alpha^4(x_7))\!+\!g(\{\alpha^2(x_3)\{x_1x_2x_4\}\alpha^2(x_5)\},\! \alpha^4(x_6),\! \alpha^4(x_7))\\
 +&g(\{\alpha^2(x_3)\alpha^2(x_4)\{x_1x_2x_5\}\},\! \alpha^4(x_6),\! \alpha^4(x_7))\!-\!g(\{\alpha^2(x_1)\alpha^2(x_2)\{x_3x_4x_5\}\},\! \alpha^4(x_6),\! \alpha^4(x_7))\\
 +&g(\{\alpha^4(x_5),\! \{\{x_1x_2x_3\}\alpha^2(x_4)\alpha^2(x_6)\},\! \alpha^4(x_7))\!+\!g(\{\alpha^4(x_5),\! \{\alpha^2(x_3)\{x_1x_2x_4\}\alpha^2(x_6)\},\! \alpha^4(x_7))\\
 +&g(\{\alpha^4(x_5),\! \{\alpha^2(x_3)\alpha^2(x_4)\{x_1x_2x_6\}\},\! \alpha^4(x_7))\!-\!g(\{\alpha^4(x_5),\! \{\alpha^2(x_1)\alpha^2(x_2)\{x_3x_4x_6\}\},\! \alpha^4(x_7))\\
 +&g(\alpha^4(x_5), \alpha^4(x_6), \{\{x_1x_2x_3\}\alpha^2(x_4)\alpha^2(x_7)\})+g(\alpha^4(x_5), \alpha^4(x_6), \{\alpha^2(x_3)\{x_1x_2x_4\}\alpha^2(x_7)\})\\
 +&g(\alpha^4(x_5), \alpha^4(x_6), \{\alpha^2(x_3)\alpha^2(x_4)\{x_1x_2x_7\}\})-g(\alpha^4(x_5), \alpha^4(x_6), \{\alpha^2(x_1)\alpha^2(x_2)\{x_3x_4x_7\}\})\\
 -&g(\alpha^4(x_3), \alpha^4(x_4), \{\{x_1x_2x_5\}\alpha^2(x_6)\alpha^2(x_7)\})-g(\alpha^4(x_3), \alpha^4(x_4), \{\alpha^2(x_5)\{x_1x_2x_6\}\alpha^2(x_7)\})\\
 -&g(\alpha^4(x_3), \alpha^4(x_4), \{\alpha^2(x_5)\alpha^2(x_6)\{x_1x_2x_7\}\})+g(\alpha^4(x_3), \alpha^4(x_4), \{\alpha^2(x_1)\alpha^2(x_2)\{x_5x_6x_7\}\})\\
 +&g(\alpha^4(x_1), \alpha^4(x_2), \{\{x_3x_4x_5\}\alpha^2(x_6)\alpha^2(x_7)\})+g(\alpha^4(x_1), \alpha^4(x_2), \{\alpha^2(x_5)\{x_3x_4x_6\}\alpha^2(x_7)\})\\
 +&g(\alpha^4(x_1), \alpha^4(x_2), \{\alpha^2(x_5)\alpha^2(x_6)\{x_3x_4x_7\}\})-g(\alpha^4(x_1), \alpha^4(x_2), \{\alpha^2(x_3)\alpha^2(x_4)\{x_5x_6x_7\}\})\\
 =&0,
 \end{align*}
where the items that could be canceled in pairs are omitted. The proof is completed.
\epf
Define
$$HomZ\!^1(L,\!L)\!\times\! HomZ\!^1(L,\!L)\!=\!\{(f,\! f)\!\in \! HomC\!^1(L,\!L)\!\times\! HomC\!^1(L,\!L)~|~(\delta_I^1,\! \delta_{II}^1)(f,\!f)=(0,\! 0)\},$$
\begin{align*}
&HomZ^2(L,L)\times HomZ^3(L,L)\\
=&\{(f,g)\in HomC^2(L,L)\times HomC^3(L,L)~|~ (\delta_I^2, \delta_{II}^2)(f,g)=(d_I^2, d_{II}^2)(f,g)=(0,0)\},
\end{align*}
\begin{align*}
&HomZ^4(L,L)\times HomZ^5(L,L)\\
=&\{(f,g)\in HomC^4(L,L)\times HomC^5(L,L)~|~ (\delta_I^3, \delta_{II}^3)(f,g)=(0,0)\},
\end{align*}
$$HomB^2(L,L)\times HomB^3(L,L)=\{(\delta_I^1, \delta_{II}^1)(f,f) ~|~ f\in HomC^1(L,L)\},$$
$$HomB^4(L,L)\!\times\! HomB^5(L,L)\!=\!\{(\delta_I^2, \delta_{II}^2)(f,g) ~|~ (f,g)\in HomC^2(L,L)\!\times\! HomC^3(L,L)\}.$$
Then by Theorem \ref{delta^2=0},
$$HomB^2(L,L)\times HomB^3(L,L)\subseteq HomZ^2(L,L)\times HomZ^3(L,L),$$
$$HomB^4(L,L)\times HomB^5(L,L)\subseteq HomZ^4(L,L)\times HomZ^5(L,L).$$ So one could define
$$HomH^1(L,L)\times HomH^1(L,L)=HomZ^1(L,L)\times HomZ^1(L,L),$$
\begin{align*}
HomH^2(L,L)\times HomH^3(L,L)=\frac{HomZ^2(L,L)\times HomZ^3(L,L)}{HomB^2(L,L)\times HomB^3(L,L)}
\end{align*}
and
\begin{align*}
HomH^4(L,L)\times HomH^5(L,L)=\frac{HomZ^4(L,L)\times HomZ^5(L,L)}{HomB^4(L,L)\times HomB^5(L,L)}
\end{align*}
as the \textbf{first, second and third cohomology groups} of $(L, [\cdot,\cdot], \{\cdot, \cdot,\cdot\}, \alpha)$, respectively.
\bdefn
A linear map $D: L\rightarrow L$ is called the $\alpha^k$-derivation of $(L, [\cdot,\cdot], \{\cdot, \cdot,\cdot\}, \alpha)$, if $D$ satisfies $D\circ \alpha=\alpha\circ D$ and
\begin{align*}
D([xy])&=[\alpha^k(x)D(y)]+[D(x)\alpha^k(y)],\\
D(\{xyz\})&=\{D(x)\alpha^k(y)\alpha^k(z)\}+\{\alpha^k(x)D(y)\alpha^k(z)\}+\{\alpha^k(x)\alpha^k(y)D(z)\},
\end{align*}
where $\alpha^k=\underbrace{\alpha\circ\cdots\circ\alpha}_{k}$ and $\alpha^0=\id_L$. It is straightforward to show that $D$ is an $\alpha^0$-derivation of $(L, [\cdot,\cdot], \{\cdot, \cdot,\cdot\}, \alpha)$ if and only if $(D,D)\in HomH^1(L,L) \times HomH^1(L,L)$. Denote by $\Der_{\alpha^k}(L)$ the set of all $\alpha^k$-derivations of $(L, [\cdot,\cdot], \{\cdot, \cdot,\cdot\}, \alpha)$.
\edefn
\bthm
Set $\Der(L)=\bigoplus_{k\geq0}\Der_{\alpha^k}(L)$. Then $\Der(L)$ is a Lie algebra.
\ethm
\bpf
It is sufficient to prove $[\Der_{\alpha^k}(L), \Der_{\alpha^s}(L)]\subseteq \Der_{\alpha^{k+s}}(L)$. Now suppose that $D\in \Der_{\alpha^k}(L)$ and $D'\in \Der_{\alpha^s}(L)$. Then
$$[D,D']\circ\alpha=D\circ D'\circ\alpha-D'\circ D\circ\alpha=\alpha\circ(D\circ D'-D'\circ D)=\alpha\circ[D,D'].$$
Note that
\begin{align*}
[D,D']([xy])=&D([\alpha^s(x)D'(y)]+[D'(x)\alpha^s(y)])-D'([\alpha^k(x)D(y)]+[D(x)\alpha^k(y)])\\
            =&[D\alpha^s(x)\alpha^kD'(y)]+[\alpha^{k+s}(x)DD'(y)]+[DD'(x)\alpha^{k+s}(y)]+[\alpha^kD'(x)D\alpha^s(y)]\\
            &\!-\![D'\alpha^k(x)\alpha^sD(y)]\!-\![\alpha^{k+s}(x)D'D(y)]\!-\![D'D(x)\alpha^{k+s}(y)]\!-\![\alpha^sD(x)D'\alpha^k(y)]\\
            =&[[D,D'](x)\alpha^{k+s}(y)]+[\alpha^{k+s}(x)[D,D'](y)],
\end{align*}
and
\begin{align*}
 &[D,D'](\{xyz\})\\
=&D(\{D'(x)\alpha^s(y)\alpha^s(z)\}+\{\alpha^s(x)D'(y)\alpha^s(z)\}+\{\alpha^s(x)\alpha^s(y)D'(z)\})\\
 &-D'(\{D(x)\alpha^k(y)\alpha^k(z)\}+\{\alpha^k(x)D(y)\alpha^k(z)\}+\{\alpha^k(x)\alpha^k(y)D(z)\})\\
=&\{DD'(x)\alpha^{k+s}(y)\alpha^{k+s}(z)\}+\{\alpha^kD'(x)D\alpha^s(y)\alpha^{k+s}(z)\}+\{\alpha^kD'(x)\alpha^{k+s}(y)D\alpha^s(z)\}\\
 &+\{D\alpha^s(x)\alpha^kD'(y)\alpha^{k+s}(z)\}+\{\alpha^{k+s}(x)DD'(y)\alpha^{k+s}(z)\}+\{\alpha^{k+s}(x)\alpha^kD'(y)D\alpha^s(z)\}\\
 &+\{D\alpha^s(x)\alpha^{k+s}(y)\alpha^kD'(z)\}+\{\alpha^{k+s}(x)D\alpha^s(y)\alpha^kD'(z)\}+\{\alpha^{k+s}(x)\alpha^{k+s}(y)DD'(z)\}\\
 &-\{D'D(x)\alpha^{k+s}(y)\alpha^{k+s}(z)\}-\{\alpha^sD(x)D'\alpha^k(y)\alpha^{k+s}(z)\}-\{\alpha^sD(x)\alpha^{k+s}(y)D'\alpha^k(z)\}\\
 &-\{D'\alpha^k(x)\alpha^sD(y)\alpha^{k+s}(z)\}-\{\alpha^{k+s}(x)D'D(y)\alpha^{k+s}(z)\}-\{\alpha^{k+s}(x)\alpha^sD(y)D'\alpha^k(z)\}\\
 &-\{D'\alpha^k(x)\alpha^{k+s}(y)\alpha^sD(z)\}-\{\alpha^{k+s}(x)D'\alpha^k(y)\alpha^sD(z)\}-\{\alpha^{k+s}(x)\alpha^{k+s}(y)D'D(z)\}\\
=&\{[D,D'](x)\alpha^{k+s}(y)\alpha^{k+s}(z)\}+\!\{\alpha^{k+s}(x)[D,D'](y)\alpha^{k+s}(z)\}+\!\{\alpha^{k+s}(x)\alpha^{k+s}(y)[D,D'](z)\}.
\end{align*}
It follows that $[D,D']\in \Der_{\alpha^{k+s}}(L)$.
\epf
\section{Deformations of a Hom-Lie-Yamaguti algebra}
Suppose that $(L, [\cdot,\cdot], \{\cdot, \cdot,\cdot\}, \alpha)$ is a HLYA over $\K$. Let $\K[[t]]$ be the ring of formal power series over $\K$ and $L[[t]]$ be the set of formal power series over $L$. Then for a $\K$-bilinear map $f:L\times L\rightarrow L$ and a $\K$-trilinear map $g:L\times L\times L\rightarrow L$, it is natural to extend them to be a $\K[[t]]$-bilinear map $f:L[[t]]\times L[[t]]\rightarrow L[[t]]$ and a $\K[[t]]$-trilinear map $g:L[[t]]\times L[[t]]\times L[[t]]\rightarrow L[[t]]$ by
$$f\left(\sum_{i\geq0}x_it^i,\sum_{j\geq0}y_jt^j\right)=\sum_{i,j\geq0}f(x_i,y_j)t^{i+j},$$
$$g\left(\sum_{i\geq0}x_it^i,\sum_{j\geq0}y_jt^j,\sum_{k\geq0}z_kt^k\right)=\sum_{i,j,k\geq0}g(x_i,y_j,z_k)t^{i+j+k}.$$
\bdefn
Suppose that $(L, [\cdot,\cdot], \{\cdot, \cdot,\cdot\}, \alpha)$ is a HLYA over $\K$. A one-parameter formal deformation of $(L, [\cdot,\cdot], \{\cdot, \cdot,\cdot\}, \alpha)$ is a pair of formal power series $(f_t, g_t)$ of the form $$f_t=[\cdot,\cdot]+\sum_{i\geq1}f_it^i, \quad g_t=\{\cdot, \cdot,\cdot\}+\sum_{i\geq1}g_it^i,$$
where each $f_i:L\times L\rightarrow L$ is a $\K$-bilinear map (extended to be $\K[[t]]$-bilinear) and each $g_i:L\times L\times L\rightarrow L$ is a $\K$-trilinear map (extended to be $\K[[t]]$-trilinear) such that $(L[[t]], f_t, g_t, \alpha)$ is a HLYA over $\K[[t]]$. Set $f_0=[\cdot,\cdot]$ and $g_0=\{\cdot, \cdot,\cdot\}$, then $f_t$ and $g_t$ can be written as $f_t=\sum_{i\geq0}f_it^i$ and $g_t=\sum_{i\geq0}g_it^i$, respectively.
\edefn

Note that $(L[[t]], f_t, g_t, \alpha)$ is required to be a HLYA. Then the following equations must be satisfied:
\begin{gather}
\alpha\circ f_t(x, y)=f_t(\alpha(x), \alpha(y)),\label{deq1}\\
\alpha\circ g_t(x, y, z)=g_t(\alpha(x), \alpha(y), \alpha(z)),\label{deq2}\\
f_t(x, x)=0,\label{deq3}\\
g_t(x,x,y)=0,\label{deq4}\\
\circlearrowleft_{x,y,z}(f_t(f_t(x, y), \alpha(z))+g_t(x, y, z))=0,\label{deq5}\\
\circlearrowleft_{x,y,z}g_t(f_t(x, y), \alpha(z), \alpha(u))=0,\label{deq6}\\
g_t(\alpha(x), \alpha(y), f_t(z, u))=f_t(g_t(x, y, z), \alpha^2(u))+f_t(\alpha^2(z), g_t(x, y, u)),\label{deq7}\\
\begin{aligned}
g_t(\alpha^2(u),\alpha^2(v),g_t(x,y,z))=&g_t(g_t(u,v,x),\alpha^2(y),\alpha^2(z))+g_t(\alpha^2(x),g_t(u,v,y),\alpha^2(z))\\
&+g_t(\alpha^2(x),\alpha^2(y),g_t(u,v,z)).
\end{aligned}\label{deq8}
\end{gather}
Equations (\ref{deq1})-(\ref{deq8}) are equivalent to
\begin{gather}
\alpha\circ f_n(x, y)=f_n(\alpha(x), \alpha(y)),\tag{\ref{deq1}$'$}\label{deq1'}\\
\alpha\circ g_n(x, y, z)=g_n(\alpha(x), \alpha(y), \alpha(z)),\tag{\ref{deq2}$'$}\label{deq2'}\\
f_n(x, x)=0,\tag{\ref{deq3}$'$}\label{deq3'}\\
g_n(x,x,y)=0,\tag{\ref{deq4}$'$}\label{deq4'}\\
\circlearrowleft_{x,y,z}\left(\sum_{i+j=n}f_i(f_j(x, y), \alpha(z))+ g_n(x, y, z)\right)=0,\tag{\ref{deq5}$'$}\label{deq5'}\\
\circlearrowleft_{x,y,z}\sum_{i+j=n}g_i(f_j(x, y), \alpha(z), \alpha(u))=0,\tag{\ref{deq6}$'$}\label{deq6'}\\
\sum_{i+j=n}g_i(\alpha(x), \alpha(y), f_j(z, u))=\sum_{i+j=n}\Big(f_i(g_j(x, y, z), \alpha^2(u))+f_i(\alpha^2(z), g_j(x, y, u))\Big),\tag{\ref{deq7}$'$}\label{deq7'}\\
\begin{aligned}
\sum_{i+j=n}g_i(\alpha^2(u),\alpha^2(v),g_j(x,y,z))&=\sum_{i+j=n}\Big(g_i(g_j(u,v,x),\alpha^2(y),\alpha^2(z))\\
+&g_i(\alpha^2(x),g_j(u,v,y),\alpha^2(z))\!+\!g_i(\alpha^2(x),\alpha^2(y),g_j(u,v,z))\Big),
\end{aligned}\tag{\ref{deq8}$'$}\label{deq8'}
\end{gather}
respectively. These equations are called the \textbf{deformation equations} of a HLYA. Equations (\ref{deq1'})-(\ref{deq4'}) imply that $(f_i, g_i)\in HomC^2(L,L)\times HomC^3(L,L)$.

Let $n=1$ in (\ref{deq5'})-(\ref{deq8'}). Then
\begin{gather*}
0=\circlearrowleft_{x,y,z}\left([f_1(x, y)\alpha(z)]+f_1([xy], \alpha(z))+g_1(x, y, z)\right),\\
0=\circlearrowleft_{x,y,z}\left(\{f_1(x, y)\alpha(z)\alpha(u)\}+g_1([xy], \alpha(z), \alpha(u))\right),\\
\begin{aligned}
0=&\{\alpha(x)\alpha(y)f_1(z, u)\}+g_1(\alpha(x), \alpha(y), [zu])-[g_1(x, y, z)\alpha^2(u)]\\
  &-f_1(\{xyz\}, \alpha^2(u))-[\alpha^2(z)g_1(x, y, u)]-f_1(\alpha^2(z), \{xyu\}),
\end{aligned}\\
\begin{aligned}
0=&\{\alpha^2(u)\alpha^2(v)g_1(x,y,z)\}-\{g_1(u, v, x)\alpha^2(y)\alpha^2(z)\}-\{\alpha^2(x)g_1(u,v,y)\alpha^2(z)\}\\
  &-\{\alpha^2(x)\alpha^2(y)g_1(u,v,z)\}+g_1(\alpha^2(u),\alpha^2(v),\{xyz\})-g_1(\{uvx\}, \alpha^2(y),\alpha^2(z))\\
  &-g_1(\alpha^2(x),\{uvy\},\alpha^2(z))-g_1(\alpha^2(x),\alpha^2(y),\{uvz\}),
\end{aligned}
\end{gather*}
which imply $(\delta_I^2, \delta_{II}^2)(f_1, g_1) =(d_I^2,d_{II}^2)(f_1, g_1)=(0,0)$, i.e.,
$$(f_1, g_1)\in HomZ^2(L,L)\times HomZ^3(L,L).$$
The first order term $(f_1, g_1)$ is called the \textbf{infinitesimal} of $(f_t, g_t)$.
\bdefn
Suppose that $(f_t, g_t)$ and $(f_t', g_t')$ are two one-parameter formal deformations of $(L, [\cdot,\cdot], \{\cdot, \cdot,\cdot\}, \alpha)$. They are called equivalent, denoted by $(f_t, g_t)\sim(f_t', g_t')$, if there is a linear isomorphism of HLYA $\Phi_t=\sum_{i\geq0}\phi_it^i: (L[[t]], f_t, g_t, \alpha)\rightarrow(L[[t]], f_t', g_t', \alpha)$ such that
\begin{gather*}
\phi_0=\id_L,\quad \Phi_t\circ \alpha=\alpha\circ\Phi_t,\\
\Phi_t\circ f_t(x, y)=f_t'(\Phi_t(x), \Phi_t(y)),\quad \Phi_t\circ g_t(x, y, z)=g_t'(\Phi_t(x), \Phi_t(y), \Phi_t(z)).
\end{gather*}
In the case $(f_1, g_1)=(f_2, g_2)=\cdots=(0,0)$, $(f_t, g_t)=(f_0, g_0)$ is called the null deformation. A one-parameter formal deformation $(f_t, g_t)$ is said to be trivial if $(f_t, g_t)\sim(f_0, g_0)$. A HLYA $(L, [\cdot,\cdot], \{\cdot, \cdot,\cdot\}, \alpha)$ is called analytically rigid if every one-parameter formal deformation $(f_t, g_t)$ is trivial.
\edefn
\bthm
Let $(f_t, g_t)$ and $(f_t', g_t')$ be equivalent one-parameter formal deformations of $(L, [\cdot,\cdot], \{\cdot, \cdot,\cdot\}, \alpha)$. Then $(f_1, g_1)$ and $(f_1', g_1')$ belong to the same cohomology class in $HomH^2(L,L)\times HomH^3(L,L)$.
\ethm
\bpf
It is sufficient to prove $(f_1-f_1', g_1-g_1')\in HomB^2(L,L)\times HomB^3(L,L)$. Suppose that $\Phi_t=\sum_{i\geq0}\phi_it^i: (L[[t]], f_t, g_t, \alpha)\rightarrow(L[[t]], f_t', g_t', \alpha)$ is an isomorphism such that
\begin{gather*}
\phi_0=\id_L,\quad \Phi_t\circ \alpha=\alpha\circ\Phi_t,\\
\Phi_t\circ f_t(x, y)=f_t'(\Phi_t(x), \Phi_t(y)),\quad \Phi_t\circ g_t(x, y, z)=g_t'(\Phi_t(x), \Phi_t(y), \Phi_t(z)).
\end{gather*}
Then $\phi_1\in HomC^1(L,L)$ and
\begin{gather*}
\sum_{i\geq0}\phi_it^i\left(\sum_{j\geq0}f_j(x, y)t^j\right)=\sum_{i\geq0}f_i'\left(\sum_{k\geq0}\phi_k(x)t^k, \sum_{l\geq0}\phi_l(y)t^l\right),\\
\sum_{i\geq0}\phi_it^i\left(\sum_{j\geq0}g_j(x, y, z)t^j\right)=\sum_{i\geq0}g_i'\left(\sum_{k\geq0}\phi_k(x)t^k, \sum_{l\geq0}\phi_l(y)t^l, \sum_{m\geq0}\phi_m(z)t^m\right).
\end{gather*}
Hence
\begin{gather*}
f_1(x, y)+\phi_1([xy])=[x\phi_1(y)]+[\phi_1(x)y]+f_1'(x, y),\\
g_1(x, y, z)+\phi_1(\{xyz\})=\{\phi_1(x)yz\}+\{x\phi_1(y)z\}+\{xy\phi_1(z)\}+g_1'(x, y, z).
\end{gather*}
Therefore, $(f_1-f_1', g_1-g_1')=(\delta_I^1, \delta_{II}^1)(\phi_1, \phi_1)\in HomB^2(L,L)\times HomB^3(L,L)$.
\epf
\bthm
Suppose that $(L, [\cdot,\cdot], \{\cdot, \cdot,\cdot\}, \alpha)$ is a HLYA. Then $(L, [\cdot,\cdot], \{\cdot, \cdot,\cdot\}, \alpha)$ is analytically rigid, if $HomH^2(L,L)\times HomH^3(L,L)=0$.
\ethm
\bpf
Let $(f_t, g_t)$ be a one-parameter formal deformation of $(L, [\cdot,\cdot], \{\cdot, \cdot,\cdot\}, \alpha)$. Suppose $f_t=f_0+\sum_{i\geq r}f_it^i$ and $g_t=g_0+\sum_{i\geq r}g_it^i$. Set $n=r$ in (\ref{deq5'})-(\ref{deq8'}). It follows that
$$(f_r, g_r)\in HomZ^2(L,L)\times HomZ^3(L,L)=HomB^2(L,L)\times HomB^3(L,L).$$
Then there exists $h_r\in  HomC^1(L,L)$ such that $(f_r, g_r)=(\delta_I^1h_r, \delta_{II}^1h_r)$.

Consider $\Phi_t=\id_L-h_rt^r$. Then $\Phi_t: L\rightarrow L$ is a linear isomorphism and $\Phi_t\circ \alpha=\alpha\circ\Phi_t$. Let
$$f_t'(x, y)=\Phi_t^{-1}f_t(\Phi_t(x), \Phi_t(y)), \quad g_t'(x, y, z)=\Phi_t^{-1}g_t(\Phi_t(x), \Phi_t(y), \Phi_t(z)).$$
Assume that $f_t'=\sum_{i\geq 0}f_i't^i$ and use the fact $\Phi_tf_t'(x, y)=f_t(\Phi_t(x), \Phi_t(y))$. Then
$$(\id_L-h_rt^r)\sum_{i\geq 0}f_i'(x, y)t^i=\left(f_0+\sum_{i\geq r}f_it^i\right)(x-h_r(x)t^r, y-h_r(y)t^r),$$
that is,
\begin{align*}
&\sum_{i\geq 0}f_i'(x, y)t^i-\sum_{i\geq 0}h_r\circ f_i'(x, y)t^{i+r}\\
=&f_0(x, y)-f_0(h_r(x), y)t^r-f_0(x, h_r(y))t^r+f_0(h_r(x), h_r(y))t^{2r}\\
 &+\sum_{i\geq r}f_i(x, y)t^i-\sum_{i\geq r}f_i(h_r(x), y)t^{i+r}-\sum_{i\geq r}f_i(x, h_r(y))t^{i+r}+\sum_{i\geq r}f_i(h_r(x), h_r(y))t^{i+2r}.
\end{align*}
So $f_0'(x, y)=f_0(x, y)=[xy]$, $f_1'(x, y)=\cdots=f_{r-1}'(x, y)=0$ and
$$f_r'(x, y)-h_r([xy])=-[h_r(x)y]-[xh_r(y)]+f_r(x, y).$$
Hence $f_r'(x, y)=-\delta_I^1h_r(x, y)+f_r(x, y)=0$, which implies $f_t'=[\cdot,\cdot]+\sum_{i\geq r+1}f_i't^i$. In the same way, we have $g_t'=\{\cdot, \cdot,\cdot\}+\sum_{i\geq r+1}g_i't^i$. It is clear that $(f_t', g_t')$ is a one-parameter formal deformation of $(L, [\cdot,\cdot], \{\cdot, \cdot,\cdot\}, \alpha)$ and $(f_t, g_t)\sim(f_t', g_t')$. By induction, one gets $(f_t, g_t)\sim(f_0, g_0)$. Therefore, $(L, [\cdot,\cdot], \{\cdot, \cdot,\cdot\}, \alpha)$ is analytically rigid.
\epf

In the deformation theory of other algebraic structures, the obstructions are in the same cohomology theory as the infinitesimal deformations but one dimension higher. But the usual approach doesn't work for HLYAs:

Let $(f_0, g_0)=([\cdot,\cdot], \{\cdot, \cdot, \cdot\})$ and $(f_1, g_1)\in HomZ^2(L,L)\times HomZ^3(L,L)$. Then $(f_0, g_0)$, $(f_1, g_1)$ satisfy the deformation equations (\ref{deq1'})-(\ref{deq8'}) for $n=1$. Set
\begin{gather*}
F(x, y, z, u)=f_1(g_1(x, y, z), \alpha^2(u))+f_1(\alpha^2(z), g_1(x, y, u))-g_1(\alpha(x), \alpha(y), f_1(z, u)),\\
\begin{aligned}
G(u, v, x, y, z)=&g_1(g_1(u, v, x), \alpha^2(y), \alpha^2(z))+g_1(\alpha^2(x), g_1(u, v, y), \alpha^2(z))\\
  &+g_1(\alpha^2(x), \alpha^2(y), g_1(u, v, z))-g_1(\alpha^2(u), \alpha^2(v), g_1(x, y, z)).
\end{aligned}
\end{gather*}
Then $(F, G)\in HomZ^4(L,L)\times HomZ^5(L,L)$. If $HomH^4(L,L)\times HomH^5(L,L)=0$, then there exists a pair $(f_2, g_2)\in HomC^2(L,L)\times HomC^3(L,L)$ such that $(\delta_I^2, \delta_{II}^2)(f_2, g_2)=(-F,-G)$. Based on one's experience in other algebras, $(f_0, g_0)$, $(f_1, g_1)$ and $(f_2, g_2)$ would satisfy the deformation equations (\ref{deq1'})-(\ref{deq8'}) for $n=2$. Note that (\ref{deq1'})-(\ref{deq4'}) clearly hold since $(f_i, g_i)\in HomC^2(L,L)\times HomC^3(L,L)$ and it is straightforward to verify $(f_0, g_0)$, $(f_1, g_1)$ and $(f_2, g_2)$ satisfying (\ref{deq7'}) and (\ref{deq8'}) for $n=2$. However, one could not prove that $(f_0, g_0)$, $(f_1, g_1)$ and $(f_2, g_2)$ satisfy (\ref{deq5'}) or (\ref{deq6'}) for $n=2$. Therefore, the obstructions of a HYLA involve other cohomology theory instead of the one carried over from Lie-Yamaguti algebras in \cite{Yamaguti1}, directly.

 {\bf ACKNOWLEDGEMENTS}\quad The authors would like to thank the referee for valuable comments and
suggestions on this article.


\begin{thebibliography}{99}
\bibitem{Ammar&Ejbehi&Makhlouf} F. Ammar, Z. Ejbehi and A. Makhlouf, \emph{Cohomology and deformations of Hom-algebras.} J. Lie Theory \textbf{21} (2011), no. 4, 813--836.
\bibitem{Ammar&Mabrouk&Makhlouf} F. Ammar, S. Mabrouk and A. Makhlouf, \emph{Representations and cohomology of $n$-ary multiplicative Hom-Nambu-Lie algebras.} J. Geom. Phys. \textbf{61} (2011), no. 10, 1898--1913.
\bibitem{Benayadi&Makhlouf} S. Benayadi and A. Makhlouf, \emph{Hom-Lie algebras with invariant nondegenerate bilinear forms.} J. Geom. Phys., \textbf{76} (2014), no. 2, 38--60.
\bibitem{Elhamdadi&Makhlouf} M. Elhamdadi and A. Makhlouf, \emph{Deformations of Hom-alternative and Hom-Malcev algebras.} Algebras Groups Geom. \textbf{28} (2011), no. 2, 117--145.
\bibitem{Flato&Gerstenhaber&Voronov} M. Flato, M. Gerstenhaber and A. Voronov, \emph{Cohomology and deformation of Leibniz pairs.} Lett. Math. Phys. \textbf{34} (1995), no. 1, 77--90.
\bibitem{Gaparayi&Nourou1} D. Gaparayi and A. Nourou Issa, \emph{A twisted generalization of Lie-Yamaguti algebras.} Int. J. Algebra \textbf{6} (2012), no. 5--8, 339--352.
\bibitem{Gaparayi&Nourou2} D. Gaparayi and A. Nourou Issa, \emph{Hom-Lie-Yamaguti structures on Hom-Leibniz algebras.} ArXiv: 1208.6038 (2012).
\bibitem{Gerstenhaber1} M. Gerstenhaber, \emph{On the deformation of rings and algebras.} Ann. of Math. (2) \textbf{79} (1964), 59--103.
\bibitem{Gerstenhaber2} M. Gerstenhaber, \emph{On the deformation of rings and algebras.II.} Ann. of Math. \textbf{84} (1966), 1--19.
\bibitem{Gerstenhaber3} M. Gerstenhaber, \emph{On the deformation of rings and algebras.III.} Ann. of Math. (2) \textbf{88} (1968), 1--34.
\bibitem{Gerstenhaber4} M. Gerstenhaber, \emph{On the deformation of rings and algebras.IV.} Ann. of Math. (2) \textbf{99} (1974), 257--276.
\bibitem{Hartwig&Larsson&Silvestrov} J. Hartwig, D. Larsson and S. Silvestrov, \emph{Deformations of Lie algebras using $\sigma$-derivations.} J. Algebra \textbf{295} (2006), no. 2, 314--361.
\bibitem{Hu} N. Hu, \emph{$Q$-Witt algebras, $q$-Lie algebras, $q$-holomorph structure and representations.} Algebra Colloq. \textbf{6} (1999), no. 1, 51--70.
\bibitem{Kodaira&Spencer} K. Kodaira and D. Spencer, \emph{On deformations of complex analytic structures. I, II.} Ann. of Math. (2) \textbf{67} (1958), 328--466.
\bibitem{Kubo&Taniguchi} F. Kubo and Y.  Taniguchi, \emph{A controlling cohomology of the deformation theory of Lie triple systems.} J. Algebra \textbf{278} (2004), no. 1, 242--250.
\bibitem{Larsson&Silvestrov} D. Larsson and S. Silvestrov, \emph{Quasi-hom-Lie algebras, central extensions and 2-cocycle-like identities.} J. Algebra \textbf{288} (2005), no. 2, 321--344.
\bibitem{Lin&Chen&Ma} J. Lin, L. Chen and Y. Ma, \emph{On the Deformation of Lie-Yamaguti algebras.} Acta Math. Sin. (Engl. Ser.) (Accepted).
\bibitem{Liu&Chen&Ma} Y. Liu, L. Chen and Y. Ma, \emph{Hom-Nijienhuis operators and T*-extensions of hom-Lie superalgebras.} Linear Algebra Appl. \textbf{439} (2013), no. 7, 2131--2144.
\bibitem{Makhlouf&Silvestrov} A. Makhlouf and S. Silvestrov, \emph{Notes on 1-parameter formal deformations of Hom-associative and Hom-Lie algebras.} Forum Math. \textbf{22} (2010), no. 4, 715--739.
\bibitem{Ma&Chen&Lin} Y. Ma, L. Chen and J. Lin, \emph{Cohomology and 1-parameter formal deformations of Hom-Lie triple systems.}  ArXiv:1309.3347 (2013).
\bibitem{Nijenhuis&Richardson} A. Nijenhuis and R. Richardson, \emph{Cohomology and deformations in graded Lie algebras.} Bull. Amer. Math. Soc. \textbf{72} (1966), 1--29.
\bibitem{Nomizu} K. Nomizu, \emph{Invariant affine connections on homogeneous spaces.} Amer. J. Math. \textbf{76} (1954), 33--65.
\bibitem{Sheng} Y. Sheng, \emph{Representations of hom-Lie algebras.}  Algebr. Represent. Theory \textbf{15} (2012), no. 6, 1081--1098.
\bibitem{Sheng&Chen} Y. Sheng and D. Chen, \emph{Hom-Lie 2-algebras.} J. Algebra \textbf{376} (2013), 174--195.
\bibitem{Yamaguti1} K. Yamaguti, \emph{On cohomology groups of general Lie triple systems.} Kumamoto J. Sci. Ser. A \textbf{8} (1967/1969), 135--146.
\bibitem{Yamaguti2} K. Yamaguti, \emph{On the Lie triple system and its generalization.} J. Sci. Hiroshima Univ. Ser. A  \textbf{21} (1957/1958),  155--160.
\bibitem{Yau1} D. Yau, \emph{Hom-algebras and homology.} J. Lie Theory \textbf{19} (2009), no. 2, 409--421.
\bibitem{Yau2} D. Yau, \emph{On $n$-ary Hom-Nambu and Hom-Nambu-Lie algebras.} J. Geom. Phys. \textbf{62} (2012), no. 2, 506--522.
\end{thebibliography}
\end{document}